\documentstyle[11pt]{article}

\begin{document}
\newcommand{\ol }{\overline}
\newcommand{\ul }{\underline }
\newcommand{\ra }{\rightarrow }
\newcommand{\lra }{\longrightarrow }
\newcommand{\ga }{\gamma }
\newcommand{\st }{\stackrel }
\newcommand{\scr }{\scriptsize }
\newcommand{\rt }{\rtimes }
\title{\Large\textbf{On Schur Multipliers
of Pairs and Triples of Groups with Topological Approach}}
\author{\textbf{Hanieh Mirebrahimi and Behrooz Mashayekhy \footnote{Correspondence: mashaf@math.um.ac.ir}} \\
Department of Pure Mathematics,\\ Center of Excellence in Analysis on Algebraic Structures,\\ Ferdowsi University of Mashhad,\\
P. O. Box 1159-91775, Mashhad, Iran.}
\date{ }
\maketitle
\begin{abstract}

In this paper, using a relation between Schur multipliers of pairs
and triples of groups, the fundamental group and homology groups of
a homotopy pushout of Eilenberg-MacLane spaces, we present among
other things some behaviors of Schur multipliers of pairs and triples with
respect to free, amalgamated free, and direct products and also
direct limits with topological approach.
\end{abstract}
\ \ \\
A. M. S. Classification (2000): 20J06; 55P20; 55U25; 57M07.
\newpage
\begin{center}
{\large\textbf{1. Introduction and Preliminaries}}
\end{center}
The Schur multiplier of a group $G$ is defined to be
$$M(G)=(R\cap F')/[R,F],$$ where $F/R$ is any free presentation of $G$. It is a well-known fact that $M(G)$
depends, up to isomorphism, only on $G$. Furthermore, it is easy
to see that $M(-)$ is a functor from the category of groups to the
category of abelian groups (see [5] for further details).

By a pair of groups $(G,N)$ we mean a group $G$ with a normal
subgroup $N$. A homomorphism of pairs $(G_1,N_1) \rightarrow
(G_2,N_2)$ is a group homomorphism $G_1 \rightarrow G_2$ that sends
$N_1$ into $N_2$. The Schur multiplier of a pair of groups $(G,N)$
which was first defined by G. Ellis [4] will be a functorial abelian
group $M(G,N)$ whose principal feature is a natural exact sequence
$$\cdots \rightarrow H_3(G)\rightarrow H_3(G/N)\rightarrow
M(G,N)\rightarrow M(G)\rightarrow$$
$$\ \ \ \ \ \ \ M(G/N)\rightarrow N/[N,G]\rightarrow G^{ab}\rightarrow
(G/N)^{ab}\rightarrow 0\ \ \ \ \ \ (1.1)$$ in which $H_3(G)$ is the
third homology group of $G$ with integer coefficient.

There are several possible definitions of the Schur multiplier of a
group and the Schur multiplier of a pair of groups. We are going to
deal with topological one that we present in this note.

First, we note that for any group $G$ one can construct functorially
a connected CW-complex $K(G)$, called Eilenberg-MacLane space, whose
fundamental group is isomorphic to $G$ which has all higher homotopy
groups trivial [8]. By considering $H_n(X)$ as the $n$th singular
homology group of a topological space $X$, with coefficients in the
group ${\bf Z}$, we recall the relation $H_n(G)\cong H_n(K(G))$, for
all $n\geq 0$, [1, Prop. 4.1].

By Hopf formula for any CW-complex $K$ with $\pi_1(K)=G$ and $F/R$
as a free presentation for $G$ we have the following isomorphism
$$\frac{H_2(K)}{h_2(\pi_2(K))}\cong \frac{R\cap F'}{[R,F]},$$
where $h_2$ is the corresponding Hurewicz map [3]. Hence a
topological definition of the Schur multiplier of a group $G$ can
be considered as the second homology group of the
Eilenberg-MacLane space $K(G)$, $H_2(K(G))$. This topological
interpretation of $M(G)$ can be extended to one for $M(G,N)$ as
follows.

For any two group extensions $1\rightarrow M\rightarrow
P\rightarrow Q\rightarrow 1$ and $1\rightarrow N\rightarrow
P\rightarrow R\rightarrow 1$ we consider the following homotopy
pushout
\[\begin{array}{ccccccc}
 \ \ \ \ \ K(P) \longrightarrow  K(\frac{P}{M})  \\
 \hspace{0.4cm}  \downarrow \st{}{} \hspace{1.2cm} \downarrow
 \st{}{}\\
  K(\frac{P}{N}) \longrightarrow X .\\
\end{array}\]
By Mayer-Vietoris sequence for pushout, we have the following exact
sequence $$\cdots \rightarrow H_3(P)\rightarrow H_3(Q)\oplus
H_3(R)\rightarrow H_3(X)\rightarrow H_2(P)\rightarrow $$
$$H_2(Q)\oplus H_2(R)\rightarrow H_2(X) \rightarrow H_1(P)\rightarrow H_1(Q)\oplus
H_1(R)\rightarrow $$ $$\ \ \ \ \ \ \ \ \ \ \ \ \ \ \ \ \ \ \ \
H_1(X)\rightarrow H_0(P)\rightarrow H_0(Q)\oplus H_0(R)\rightarrow
H_0(X)\rightarrow 0.\ \ \ \ \ \ \ \ \ \ (1.2)$$ Using [2, Corollary
3.4] we have $$\pi_1(X)\cong \frac{P}{MN}\ \ \ \ {\rm and}\ \ \ \
\pi_2(X)\cong \frac{M\cap N}{[M,N]}.$$ If we make the assumption
$P=MN$, then $X$ is $1-$connected and so by Hurewicz Theorem we have
$$H_1(X)=0\ \ \ \ {\rm and}\ \ \ \ H_2(X)\cong\pi_2(X)\cong \frac{M\cap
N}{[M,N]}.$$ Hence the exact sequence $(1.2)$ becomes as follows
$$\cdots \rightarrow H_3(P)\rightarrow H_3(Q)\oplus H_3(R)\rightarrow H_3(X)\rightarrow H_2(P)\rightarrow $$
$$H_2(Q)\oplus H_2(R)\rightarrow \frac{M\cap N}{[M,N]} \rightarrow
P^{ab}\rightarrow Q^{ab}\oplus R^{ab}\rightarrow 0.\ \ \ \ \ \ \ \ \
(1.3)$$ Now, in special case, if we consider the two group
extensions $1\rightarrow N\rightarrow G\rightarrow G/N\rightarrow 1$
and $1\rightarrow G\rightarrow G\rightarrow 1\rightarrow 1$
corresponding to a pair of groups $(G,N)$ and the following homotopy
pushout
\[\begin{array}{ccccccc}
 \ \ \ \ \ K(G) \longrightarrow  K(\frac{G}{N})  \\
 \hspace{0.4cm}  \downarrow \st{}{} \hspace{1cm} \downarrow
 \st{}{}\\
 \ \ \ \ \ \ \ \ \ \ \ \ \ \ \ \ \ \ \ \ 1 \longrightarrow X,\ \ \ \ \ \ \ \ \ \ (1.4) \\
\end{array}\]
then we have the following natural exact sequence as $(1.1)$
$$H_3(G)\rightarrow H_3(G/N)\rightarrow
H_3(X)\rightarrow H_2(G)=M(G)\rightarrow$$
$$H_2(G/N)=M(G/N)\rightarrow H_2(X)=N/[N,G]\rightarrow G^{ab}\rightarrow
(G/N)^{ab}\rightarrow 0.$$ Hence a topological interpretation of the
 Schur multiplier of a pair of groups $(G,N)$ can be considered as
the third homology group of a space $X$ which is the homotopy
pushout corresponding to the
pair of groups $(G,N)$ as $(1.4)$. \\
\textbf{\bf Remark 1.1.} As we mentioned before, the notion of the
Schur multiplier of a pair of groups was introduced by G. Ellis in
[4]. He presented several possible definitions of the notion through
them the topological one is as follows.

For a pair of groups $(G,N)$, the natural epimorphism $G\rightarrow
G/N$ induces functorially the continuous map $f: K(G)\rightarrow
K(G/N)$. Suppose that $M(f)$ is the mapping cylinder of $f$
containing $K(G)$ as a subspace and is also homotopically equivalent
to the space $K(G/N)$. Take $K(G,N)$ to be the mapping cone of the
cofibration $K(G)\hookrightarrow M(f)$. By Mayer-Vietoris, the
cofibration sequence $K(G)\hookrightarrow M(f)\rightarrow K(G,N)$
induces a natural long exact homology sequence $$\cdots \rightarrow
H_{n+1}(G/N)\rightarrow H_{n+1}(K(G,N))\rightarrow
H_{n}(G)\rightarrow H_{n}(G/N)\rightarrow \cdots$$ for $n\geq 0$. G.
Ellis [4] showed that the Schur multiplier of the pair $(G,N)$ can
be considered as the third homology group of the cofiber space
$K(G,N)$. We note that the mapping cone $K(G,N)$ of the cofiber
$K(G)\rightarrow M(f)$ is homotopically equivalent to the space $X$
which is the homotopy pushout corresponding to the pair of groups
$(G,N)$. Therefore our topological interpretation of the Schur
multiplier of a pair of groups $(G,N)$ is equivalent to the
topological definition of G. Ellis.\\

In section two, using the topological interpretations, we present
among other things some behaviors of the Schur multiplier of pairs
of the free product, the amalgamated free product, and the direct
product. Also, we show that the Schur multiplier of a pair commutes
with direct limits in some cases.

By a triple of groups $(G,M,N)$ we mean a group $G$ with two normal
subgroups $M$ and $N$. A homomorphism of triples
$(G_1,M_1,N_1)\rightarrow (G_2,M_2,N_2)$ is a group homomorphism
$G_1\rightarrow G_2$ that sends $M_1$ into $M_2$ and $N_1$ into
$N_2$. G. Ellis $[4]$ defined the Schur multiplier of a triple
$(G,M,N)$ as a functorial abelian group $M(G,M,N)$ whose principle
feature is a natural exact sequence
$$\cdots\rightarrow
H_3(G,N)\rightarrow H_3(G/M,MN/M)\rightarrow M(G,M,N)\rightarrow
M(G,N)$$$$M(G/M,MN/M)\rightarrow M\cap N/[M\cap N,G][M,N]$$$$\ \ \ \
\ \ \ \ \ \ \ \ \ \rightarrow N/[N,G]\rightarrow
NM/M[N,G]\rightarrow 0.\ \ \ \ \ \ \ \ \ \ \ \ (1.5)$$ He also gave
a topological interpretation for $M(G,M,N)$. In section three, first
we give a topological definition for the Schur multiplier of a
triple which is equivalent to the one of Ellis. Then we show that
the Schur multiplier of a triple commutes with direct limits with
some conditions. Second, we define a new version of the Schur
multiplier of a triple $(G,M,N)$ which is more natural
generalization of the Schur multiplier of a pair $(G,N)$ than the
one of Ellis. We show that our new notion is coincide with the one
of Ellis if $G=MN$. Finally, we present behaviors of this new
version of the Schur multiplier of a triple with respect to free,
amalgamated free and direct products and also a better behavior with
respect to direct limits than the one of Ellis.

\begin{center} {\large\textbf{2. Schur Multipliers of Pairs}}
\end{center}
The following two known results obtained from the existence of the
natural long exact sequence $(1.1)$. Using the properties of
Eilenberg-MacLane spaces, we reprove them by topological viewpoints.\\
\textbf{\bf Theorem 2.1.} \textit{The Schur multiplier of a group
$G$ is a special case of the Schur multiplier of a pair of groups
that is
$M(G,G)\cong M(G)$.}\\
\textbf{ Proof.} First, we know that Eilenberg-MacLane space
$K(G/G)$ can be considered contractible. Hence using the fact that
$H_i(K(G/G))=1$, for $i\geq 1$, and also by the Mayer-Vietoris
sequence due to the corresponding pushout diagram, the result holds. $\Box$\\
\textbf{\bf Theorem 2.2.} \textit{For any group $G$, the Schur
multiplier of the pair $(G,1)$ is trivial.}\\
\textbf{ Proof.} First, we recall that $M(G,1)$ is equal to the
third homology group $H_3(X)$, where $X$ is the mapping cone of the
cofibration $i:K(G)\hookrightarrow K(G)$. Hence $X$ is the quotient
space of $(K(G)\times I)\sqcup K(G)$ with respect to the equivalent
relation $(x,0)\sim (x',0)$ , $(x,1)\sim i(x)$, for all $x,x'\in
K(G)$. Thus the space $X$ is a contractible and so with trivial
homology groups which completes the proof. $\Box$

We recall that by Miller Formula $[5]$ for any two groups $G_1$
and $G_2$ we have $M(G_1*G_2)\cong M(G_1)\oplus M(G_1)$. The authors have also
presented a topological proof for this note in $[6]$. In the following, we prove a new result for the structure of the
Schur multiplier of pairs of the free product with topological method.\\
\textbf{\bf Theorem 2.3.} \textit{For any two groups $G_1$, $G_2$
and their normal subgroups $N_i\unlhd G_i$ ($i=1,2$) we have the
following isomorphism}
$$M(G_1*G_2,
 \langle N_1*N_2\rangle^{G_1*G_2})\cong M(G_1,N_1)\oplus M(G_2,N_2),$$ where
 $\langle N_1*N_2\rangle^{G_1*G_2}$ is the normal closure of $N_1*N_2$ in
 $G_1*G_2$.\\
\textbf{ Proof.} Suppose that $X_i$ $(i=1,2)$ is the pushout of
corresponding diagram
\[\begin{array}{ccccccc}
 \ \ \ \ \ K(G_i) \longrightarrow  K(\frac{G_i}{N_i})  \\
 \hspace{0.4cm}  \downarrow \st{}{} \hspace{1cm} \downarrow
 \st{}{}\\
 \ \ \ \ 1 \longrightarrow X_i. \\
\end{array}\]
Using the fact that any two direct limits commute, we conclude that
the space $X_1\vee X_2$ is also a pushout for the
following diagram
\[\begin{array}{ccccccc}
 \ \ \ \ \ \ K(G_1)\vee K(G_2) \longrightarrow  K(\frac{G_1}{N_1})\vee K(\frac{G_2}{N_2})  \\
 \hspace{0.6cm}  \downarrow \st{}{} \hspace{1.3cm} \downarrow
 \st{}{}\\
\ \ \ \ \ \ \ \ \ \ \ 1 \longrightarrow X_1\vee X_2. \\
\end{array}\]
According to Van-Kampen Theorem, we rewrite the above diagram as
follows
\[\begin{array}{ccccccc}
 \ \ \ \ \ \ \ \ \ \ \ \ \ \ \ K(G_1*G_2) \longrightarrow  K(\frac{G_1}{N_1}*\frac{G_2}{N_2})=K(\frac{G_1*G_2}{\langle N_1*N_2\rangle^{G_1*G_2}})  \\
 \hspace{0cm}  \downarrow \st{}{} \hspace{1.2cm} \downarrow
 \st{}{}\\
 \ \ \ \ \ \ 1 \longrightarrow X_1\vee X_2. \\
\end{array}\]
Now by the definition and the Mayer-Vietoris sequence for the above
diagram we have
$$M(G_1*G_2,
 \langle N_1*N_2\rangle^{G_1*G_2})=H_3(X_1\vee X_2)$$$$\cong H_3(X_1)\oplus H_3(X_2)=M(G_1,N_1)\oplus
M(G_2,N_2).\ \Box$$

We recall that the authors [6] proved that if $G$ is the
free amalgamated product of its two subgroups $G_1$ and $G_2$ over a
subgroup $H$, then the following exact sequence holds
$$\cdots \rightarrow M(H)\rightarrow M(G_1)\oplus M(G_2)\rightarrow M(G)
$$ $$\rightarrow H_{ab} \rightarrow G_{1 ab}\oplus G_{2
ab}\rightarrow G_{ab}\rightarrow \cdots .$$ Now we give the following
result for the Schur multiplier of a pair of amalgamated free product of groups by
topological argument.\\
 \textbf{\bf Theorem 2.4.} \textit{Let $(G_1,N_1)$,
$(G_2,N_2)$ be two pairs of groups and $H\leq N_1\cap N_2$. Then the
Schur multiplier of the pair of amalgamated free products
$(G_1*_H G_2 , N_1*_H N_2)$, satisfies the following exact
sequence}
$$\cdots\rightarrow H_3(K(G_1* G_2))\rightarrow H_3(\frac{G_1}{N_1})\oplus
H_3(\frac{G_2}{N_2})\rightarrow $$ $$M(G_1*_H G_2,N_1*_H
N_2)\rightarrow M(G_1*_H G_2)\rightarrow
M(\frac{G_1}{N_1})\oplus M(\frac{G_2}{N_2})\rightarrow\cdots .$$
\textbf{ Proof.} In order to introduce the Eilenberg-MacLane space
corresponding to the group $G_1*_H G_2$, first consider the space $K(H)$ and
then by attaching cells to this space in suitable ways, construct
spaces $K(G_1)$ and $K(G_2)$ (for further details see [6, Theorem2.5]). In this case, using Van-Kampen
Theorem, we have the isomorphism $\pi_1(K(G_1)\cup K(G_2))\cong
G_1*_H G_2$, and by the way of constructing the spaces $K(G_1)$
and $K(G_2)$ we can consider $K(G_1)\cup K(G_2)$ as an
Eilenberg-MacLane space
for the group $G_1*_H G_2$.\\
Also we recall that for any $i\in\{1,2\}$, $M(G_i,N_i)$ is the third
homology group $H_3(X_i)$, where $X_i$ is pushout of the following
diagram
\[\begin{array}{ccccccc}
 \ \ \ \ \ K(G_i) \longrightarrow  K(\frac{G_i}{N_i})  \\
 \hspace{0.4cm}  \downarrow \st{}{} \hspace{1cm} \downarrow
 \st{}{}\\
 \ \ \ \ \ 1 \longrightarrow X_i. \\
\end{array}\]
Using the fact that union preserves direct limit of spaces, the
following diagram is also a pushout diagram
\[\begin{array}{ccccccc}
 \ \ \ \ \ K(G_1)\cup K(G_2) \longrightarrow  K(\frac{G_1}{N_1})\cup K(\frac{G_2}{N_2})  \\
 \hspace{0.8cm}  \downarrow \st{}{} \hspace{1.5cm} \downarrow
 \st{}{}\\
 \ \ \ \ \ \ \ \ \ \ \ \ 1 \longrightarrow X_1\cup X_2. \\
\end{array}\]
Note that by the assumption $H\subseteq N_1\cap N_2$,
$K(\frac{G_1}{N_1})$ and $K(\frac{G_2}{N_2})$ has one-point
intersection. Roughly speaking, $K(\frac{G_1}{N_1})\cup
K(\frac{G_2}{N_2})$ is indeed a wedge space and so it is an
Eilenberg-MacLane space corresponding to the free product
$\frac{G_1}{N_1}*\frac{G_2}{N_2}$. On the other hand, by the
following isomorphism
$$\frac{G_1*_H G_2}{N_1*_H N_2}\cong
\frac{G_1}{N_1}*\frac{G_2}{N_2},$$ the two Eilenberg-MacLane spaces
$K(\frac{G_1*_H G_2}{N_1*_H N_2})$ and
$K(\frac{G_1}{N_1}*\frac{G_2}{N_2})=K(\frac{G_1}{N_1})\cup
K(\frac{G_2}{N_2})$ are homotopic. Hence we can rewrite above
diagram as follows
\[\begin{array}{ccccccc}
\ \ \ \ \ K(G_1*_H G_2) \longrightarrow  K(\frac{G_1*_H G_2}{N_1*_H N_2})  \\
 \hspace{0.7cm}  \downarrow \st{}{} \hspace{1.3cm} \downarrow
 \st{}{}\\
 \ \ \ \ \ \ \ \ \ \ \ \ 1 \longrightarrow X_1\cup X_2. \\
\end{array}\]
Now using the Mayer-Vietoris sequence for this recent pushout
diagram, we obtain the following exact sequence
$$\cdots\rightarrow H_3(K(G_1*_H G_2))\rightarrow H_3(K(\frac{G_1*_H
G_2}{N_1*_H N_2}))\rightarrow H_3(X_1\cup X_2)\rightarrow
$$$$H_2(K(G_1*_H G_2))\rightarrow H_2(K(\frac{G_1*_H
G_2}{N_1*_H N_2}))\rightarrow \cdots.$$  Also using
Mayer-Vietoris sequence for join spaces, we have the isomorphism
$H_n(K(\frac{G_1}{N_1})\cup K(\frac{G_2}{N_2}))\cong
H_n(K(\frac{G_1}{N_1}))\oplus H_n(K(\frac{G_2}{N_2}))$, for any
$n\in {\bf N}$. Finally, by the topological definitions of the Schur
multiplier of a group and a pair of groups, $M(G_1*_H
G_2,N_1*_H N_2)\cong H_3(X_1\cup X_2)$ and $M(G_1*_H
G_2)\cong H_2(K(G_1*_H G_2))$. Hence we get the following exact
sequence
$$\cdots\rightarrow H_3(K(G_1*_H G_2))\rightarrow H_3(\frac{G_1}{N_1})\oplus
H_3(\frac{G_2}{N_2})\rightarrow $$$$M(G_1*_H G_2,N_1*_H
N_2)\rightarrow M(G_1*_H G_2)\rightarrow
M(\frac{G_1}{N_1})\oplus M(\frac{G_2}{N_2})\rightarrow\cdots.\ \ \
\Box$$

Note that the authors [6], using topological methods, proved that for any two groups $G_1$ and $G_2$, the
following isomorphism holds
$$M(G_1\times G_2)\cong M(G_1)\oplus M(G_2)\oplus
(G_1)_{ab}\otimes (G_2)_{ab}.$$  In the
following we extend this result.\\
\textbf{\bf Theorem 2.5.} \textit{For any two pairs of groups
$(G_1,N_1)$, $(G_2,N_2)$, the Schur multiplier of a pair of the
direct products $(G_1\times G_2,N_1\times N_2)$ satisfies the
following exact sequence}
$$\cdots\rightarrow H_3(G_1)\oplus(M(G_1)\otimes G_2^{ab})\oplus
(G_1^{ab}\otimes M(G_2))\oplus $$$$H_3(G_2)\oplus
Tor(G_1^{ab},G_2^{ab})\rightarrow
H_3(\frac{G_1}{N_1})\oplus(M(\frac{G_1}{N_1})\otimes
\frac{G_2^{ab}N_2}{N_2})\oplus $$$$(\frac{G_1^{ab}N_1}{N_1}\otimes
M(\frac{G_2}{N_2}))\oplus H_3(\frac{G_2}{N_2})\oplus
 Tor(\frac{G_1^{ab}N_1}{N_1},\frac{G_2^{ab}N_2}{N_2}) \rightarrow$$$$ M(G_1\times G_2,N_1\times N_2)\rightarrow
M(G_1\times G_2)\rightarrow
M(\frac{G_1}{N_1}\times\frac{G_2}{N_2})\rightarrow\cdots .$$
\textbf{ Proof.} First, we consider the following pushout diagram
\[\begin{array}{ccccccc}
 \ \ \ \ \ K(G_1\times G_2) \longrightarrow  K(\frac{G_1\times G_2}{N_1\times N_2})  \\
 \hspace{0.4cm}  \downarrow \st{}{} \hspace{1cm} \downarrow
 \st{}{}\\
 \ \ \ \ 1 \longrightarrow X. \\
\end{array}\]
By the definition, we know that $M(G_1\times G_2,N_1\times
N_2)=H_3(X)$. On the other hand, by Mayer-Vietoris sequence for the
above diagram we conclude the following exact sequence
$$\cdots \rightarrow H_3(K(G_1\times G_2))\rightarrow H_3(K(\frac{G_1\times G_2}{N_1\times N_2}))
\rightarrow H_3(X)\rightarrow $$$$H_2(K(G_1\times G_2))\rightarrow
H_2(K(\frac{G_1\times G_2}{N_1\times N_2}))\rightarrow\cdots.$$
Using $K(G_1\times G_2)=K(G_1)\times K(G_2)$, the K$\ddot{u}$nneth
Formula and some properties of the functor $Tor$ and tensor product
we have
$$ H_3(K(G_1\times G_2))=H_3(K(G_1)\times K(G_2))\cong$$$$
H_3(K(G_1))\oplus (H_2(K(G_1))\otimes H_1(K(G_2)))\oplus
(H_1(K(G_1))\otimes H_2(K(G_2)))\oplus$$$$H_3(K(G_2))\oplus
Tor(H_1(K(G_1)),H_1(K(G_2))).$$ By the similar argument for
$H_3(K(\frac{G_1\times G_2}{N_1\times
N_2}))=H_3(K(\frac{G_1}{N_1}\times \frac{G_2}{N_2}))$ and the
isomorphisms $H_1(K(G))\cong G^{ab}$, $H_2(K(G))\cong M(G)$ we obtain the
following exact sequence $$\cdots\rightarrow
H_3(K(G_1))\oplus(M(G_1)\otimes G_2^{ab})\oplus (G_1^{ab}\otimes
M(G_2))\oplus $$$$H_3(K(G_2))\oplus
Tor(G_1^{ab},G_2^{ab})\rightarrow
H_3(K(\frac{G_1}{N_1}))\oplus(M(\frac{G_1}{N_1})\otimes
\frac{G_2^{ab}N_2}{N_2})\oplus $$$$(\frac{G_1^{ab}N_1}{N_1}\otimes
M(\frac{G_2}{N_2}))\oplus H_3(K(\frac{G_2}{N_2}))\oplus
 Tor(\frac{G_1^{ab}N_1}{N_1},\frac{G_2^{ab}N_2}{N_2}) \rightarrow$$$$ M(G_1\times G_2,N_1\times N_2)\rightarrow
M(G_1\times G_2)\rightarrow
M(\frac{G_1}{N_1}\times\frac{G_2}{N_2})\rightarrow\cdots .\Box$$
\textbf{\bf Remarks 2.6.} Using the isomorphism $H_3({\bf Z}_m)\cong {\bf Z}_m$ [1],
we can rewrite the above  exact sequence
in some special cases as follows:\\
$(i)$ If $\frac{G_i}{N_i}$ is a finite cyclic group of order $m_i$
($i=1,2$), then we have the following exact sequence $$ \cdots
\rightarrow {\bf Z}_{m_1}\oplus {\bf Z}_{m_2}\oplus {\bf
Z}_{d}\rightarrow M(G_1\times G_2,N_1\times N_2)\rightarrow $$ $$
M(G_1)\oplus M(G_2)\oplus (G_1\otimes G_2)\rightarrow {\bf
Z}_{d}\rightarrow\cdots , $$  where $d=gcd(m_1,m_2)$.
 Moreover, if
$G_1$ and $G_2$ are also cyclic, then we get the following exact
sequence $$ \cdots \rightarrow {\bf Z}_{m_1}\oplus {\bf
Z}_{m_2}\oplus {\bf Z}_{d}\rightarrow M(G_1\times G_2,N_1\times
N_2)\rightarrow $$$$ G_1\otimes G_2\rightarrow {\bf
Z}_{d}\rightarrow\cdots .$$
$(ii)$ If
$gcd(|\frac{G_1}{N_1}|,|\frac{G_2}{N_2}|)=1$, then the following
exact sequence holds $$ \cdots \rightarrow
H_3(\frac{G_1}{N_1})\oplus H_3(\frac{G_2}{N_2}) \rightarrow
M(G_1\times G_2,N_1\times N_2)\rightarrow $$$$ M(G_1)\oplus
M(G_2)\oplus (G_1\otimes G_2)\rightarrow M(\frac{G_1}{N_1})\oplus
M(\frac{G_2}{N_2})\rightarrow\cdots .$$
Moreover, if
$gcd(|G_1|,|G_2|)=1$, then we conclude the following exact sequence
$$ \cdots \rightarrow H_3(\frac{G_1}{N_1})\oplus H_3(\frac{G_2}{N_2})
\rightarrow M(G_1\times G_2,N_1\times N_2)\rightarrow $$$$
M(G_1)\oplus M(G_2)\rightarrow M(\frac{G_1}{N_1})\oplus
M(\frac{G_2}{N_2})\rightarrow\cdots .$$
$(iii)$ Finally, if $G_i$ is a finite cyclic group ($i=1,2$) such that
$gcd(|G_1|,|G_2|)=1$, then we have the following exact sequence
$$ \cdots \rightarrow {\bf Z}_{m_1}\oplus {\bf Z}_{m_2}
\rightarrow M(G_1\times G_2,N_1\times N_2)\rightarrow 0.$$

The following theorem has been proved by Ellis $[4]$. Here we
present a
topological proof for this note.\\
\textbf{\bf Theorem 2.7.} \textit{If $G=N\!>\!\!\!\!\lhd Q$ is the
semidirect product of a normal subgroup $N$ by a subgroup $Q$, then
 }\\ $$M(G)\cong M(G,N)\oplus M(Q).$$
\textbf{ Proof.} Using the corresponding pushout diagram
\[\begin{array}{ccccccc}
 \ \ \ \ \ \ \ \ \ \ \ \ \  K(G) \longrightarrow  K(\frac{G}{N})=K(Q)  \\
 \hspace{0cm}  \downarrow \st{}{} \hspace{1cm} \downarrow
 \st{}{}\\
  \ 1 \longrightarrow X, \\
\end{array}\]
and the Mayer-Vietoris sequence for this diagram, we conclude the
following exact sequence $$\cdots\rightarrow H_3(K(G))\rightarrow
H_3(K(Q))\rightarrow H_3(X)\rightarrow
$$$$ \ \
\ \ \ \ \ \ \ \ \ \ H_2(K(G))\rightarrow H_2(K(Q))\rightarrow
H_2(X)\rightarrow\cdots .\ \ \ (2.1)$$ Since $G=N\!>\!\!\!\!\lhd Q$,
there exists a homomorphism $\beta:Q\rightarrow G$ such that $\alpha
\circ \beta =1_Q$. $H_n(K(-))$ is the composition of two functors
and so the induced homomorphism $\alpha_*:H_n(K(G))\rightarrow
H_n(K(Q))$ is surjection; and hence the exactness of $(2.1)$ implies
the injectivity of the homomorphism $H_n(X)\rightarrow
H_{n-1}(K(G))$. Thus we have the following exact split sequence
$$0\rightarrow H_3(X)\rightarrow H_{2}(K(G))\rightarrow
H_{2}(K(Q))\rightarrow 0$$ which completes the proof. $\Box$\\
\textbf{\bf Theorem 2.8.} \textit{Suppose that $M$ and $N$ are two
subgroups of a group $G$ so that $M\cong MN$, then there exists the
following isomorphism }\\ $$M(MN,N)\cong M(M,M\cap N).$$
\textbf{Proof.} By the second Isomorphism Theorem we have $\frac{MN}{N}\cong
\frac{M}{M\cap N}$. Because of the functorial property of $K(-)$, we
conclude the homotopy equivalences $K(M)\approx K(MN)$ and
$K(\frac{MN}{N})\approx K(\frac{M}{M\cap N})$. Therefore, by the
uniqueness of the pushout in the category \textit{hTop}, we can
identify two following homotopoy pushout diagrams in this category
\[\begin{array}{ccccccc}
 \ \ \ \ \ K(MN) \longrightarrow  K(\frac{MN}{N})  \ \ \ \ \ K(M) \longrightarrow  K(\frac{M}{M\cap N})\\
 \hspace{0.5cm}  \downarrow \st{}{} \hspace{1cm} \downarrow \hspace{2cm}\ \ \ \ \downarrow \st{}{} \hspace{1cm} \downarrow
 \st{}{} \st{}{}\\
  \ \ \ \ \ 1 \longrightarrow X, \hspace{2cm}\ \ \  1 \longrightarrow Y.
\end{array}\]
Hence by the definition, $$M(MN,N)=H_3(X)\cong H_3(Y)=M(M,M\cap N).\
\ \ \Box$$

Note that applying the five lemma and the properties of the direct
limit of a directed system to an exact sequence, namely the fact
that this functor is exact and it commutes with homology, one can
show that for a directed system $\{(G_i,N_i)\}_{i\in I}$ of any pair
of groups, the following isomorphism holds
$$\displaystyle{M(\lim_{\rightarrow}G_i,\lim_{\rightarrow}N_i)\cong\lim_{\rightarrow}M(G_i,N_i)}.$$
However in the following we establish a topological proof for this
fact in special case, where $\{(G_i,N_i)\}_{i\in I}$ is a directed
system of pairs of abelian groups. First we need the
following lemma.\\
\textbf{\bf Lemma 2.9.} \textit{The direct limit of a direct system
$\{X_i\}_{i\in I}$ of Eilenberg-MacLane spaces is an
Eilenberg-MacLane space. Moreover, if the system $\{X_i\}_{i\in I}$
is directed and $\pi_1(X_i)$ is abelian group, for any $i\in I$,
then the Eilenberg-MacLane space $\displaystyle{\lim_{\rightarrow}
X_i}$ is corresponding to the group
$\displaystyle{\lim_{\rightarrow}\pi_1(X_i)}$.}\\
\textbf{ Proof.} Consider the induced direct system
$\{\pi_n(X_i)\}_{i\in I}$. In order to prove that
$\displaystyle{\lim_{\rightarrow} X_i}$ is an Eilenberg-MacLane space, first we
show that $\displaystyle{\pi_1(\lim_{\rightarrow} X_i)}$ is abelian. For any
$\alpha$ and $\beta$ in $\displaystyle{\pi_1(\lim_{\rightarrow} X_i)}$,
there exists two paths $\gamma$ and $\lambda$
 in $\displaystyle{\lim_{\rightarrow} X_i}$ such that $\alpha=[\gamma]$ and $\beta=[\lambda]$.
 Note that $\displaystyle{\lim_{\rightarrow}X_i}$ is a quotient space of the wedge space $\vee_{i\in I} X_i$, and
so we can consider $\gamma\cap X_i$ and $\lambda\cap X_i$ as two
paths in $X_i$, for any $i\in I$. Since $\pi_1(X_i)$'s are abelian
groups, $\gamma\cap X_i$ and $\lambda\cap X_i$ commute with each other in
the space $X_i$ $(i\in I)$ up to homotopy. Therefore two paths $\gamma$ and
$\lambda$ commute with each other in the space $\displaystyle{\lim_{\rightarrow} X_i}$ up to homotopy.
Thus $\displaystyle{\pi_1(\lim_{\rightarrow} X_i)}$ is an abelian group. By a similar
argument we can show that $\displaystyle{\pi_n(\lim_{\rightarrow} X_i)}$
is a trivial group, for any $n\geq 2$. If
$\alpha$ is an element of $\displaystyle{\pi_n(\lim_{\rightarrow} X_i)}$,
there exists an $n$-loop $\gamma$
 in $\displaystyle{\lim_{\rightarrow} X_i}$ such that $\alpha=[\gamma]$. Similar to the
 previous note, $\gamma\cap X_i$ is an $n$-loop in $X_i$ $(i\in I)$;
 and since $\pi_n(X_i)$'s are trivial groups, so $[\gamma\cap X_i]$ is trivial in
 $\pi_n(X_i)$, for any $i\in I$. Because of the structure of the space
 $\displaystyle{\lim_{\rightarrow} X_i}$, $[\gamma]$ is also trivial in
 $\displaystyle{\pi_n(\lim_{\rightarrow} X_i)}$.
 This note is true for any arbitrary $n\geq 2$ and so
  $\displaystyle{\pi_n(\lim_{\rightarrow} X_i)}$ ($n\geq 2$) is trivial.
  So if the space $X_i$, for any $i\in I$, is Eilenberg-MacLane, then the
space $ \displaystyle{\lim_{\rightarrow} X_i}$ should be also
Eilenberg-MacLane.

Moreover, if $\pi_1(X_i)$'s are abelian groups, we have
$\pi_1(X_i)\cong H_1(X_i)$, for any $i\in I$. Hence using the fact
that homology functors commute with direct limits of directed
systems, we conclude that
$$\lim_{\rightarrow}\pi_1(X_i)\cong \lim_{\rightarrow}H_1(X_i)\cong
H_1(\lim_{\rightarrow}X_i).$$ Also,
$\displaystyle{\pi_1(\lim_{\rightarrow}} X_i)$, as a homomorphic
image of the abelian group
$\displaystyle{\lim_{\rightarrow}\pi_1(X_i)}$, is abelian, and so we
have
$$H_1(\lim_{\rightarrow}X_i)\cong \pi_1(\lim_{\rightarrow}X_i)$$
which completes the proof. $\Box$\\
\textbf{\bf Theorem 2.10.} \textit{Let $\{(G_i,N_i)\}_{i\in I}$ be a
given directed system of pairs of abelian groups, then
$\displaystyle{M(\lim_{\rightarrow}G_i,\lim_{\rightarrow}N_i)\cong
\lim_{\rightarrow}M(G_i,N_i)}$.}\\
\textbf{ Proof.} First, for any $i\in I$, we consider the
corresponding pushout diagram
\[\begin{array}{ccccccc}
 \ \ \ \ \ K(G_i) \longrightarrow  K(\frac{G_i}{N_i})  \\
 \hspace{0.4cm}  \downarrow \st{}{} \hspace{1cm} \downarrow
 \st{}{}\\
 \ \ \ \ 1 \longrightarrow X_i. \\
\end{array}\]
Using the fact that any two direct limits commute, we conclude that
the following diagram is also a pushout diagram,
\[\begin{array}{ccccccc}
 \ \ \ \ \ \displaystyle{\lim_{\rightarrow} K(G_i) \longrightarrow  \lim_{\rightarrow}K(\frac{G_i}{N_i}) } \\
 \hspace{0.4cm}  \downarrow \st{}{} \hspace{1.2cm} \downarrow
 \st{}{}\\
 \ \ \ \ \ \ \ 1 \longrightarrow \displaystyle{\lim_{\rightarrow}X_i.} \\
\end{array}\]
Now by Lemma $2.9$, we rewrite the above diagram as follow
\[\begin{array}{ccccccc}
 \ \ \ \ \ \displaystyle{K(\lim_{\rightarrow}G_i) \longrightarrow  K(\lim_{\rightarrow}\frac{G_i}{N_i})}  \\
 \hspace{0.4cm}  \downarrow \st{}{} \hspace{1.2cm} \downarrow
 \st{}{}\\
 \ \ \ \ \ \ \ 1 \longrightarrow \displaystyle{\lim_{\rightarrow}X_i.} \\
\end{array}\]
Finally, by the definition and the fact that homology groups commute
with the direct limit of a directed system, we have
$$M(\lim_{\rightarrow}G_i,\lim_{\rightarrow}N_i)=
H_3(\lim_{\rightarrow}X_i)$$$$\cong \lim_{\rightarrow}H_3(X_i)=
\lim_{\rightarrow}M(G_i,N_i).\ \Box$$

\begin{center}
{\large\textbf{3. Schur Multipliers of Triples }}
\end{center}
Let $(G,M,N)$ be a triple of groups. Consider the following homotopy
pushout
\[\begin{array}{ccccccc}
 \ \ \ \ \ \ K(G,N) \longrightarrow  K(\frac{G}{M},\frac{MN}{M})  \\
 \hspace{0.4cm}  \downarrow \st{}{} \hspace{1cm} \downarrow
 \st{}{}\\
 \ \ \ 1 \longrightarrow X.\\
\end{array}\]
Using the Mayer-Vietoris sequence for homotopy pushout, we have the
following exact sequence
$$H_4(K(G,N))\rightarrow H_4(K(G/M,MN/M))\rightarrow
H_4(X)\rightarrow
$$$$H_3(K(G,N))\rightarrow H_3(K(G/M,MN/M))\rightarrow
$$$$H_3(X)\rightarrow H_2(K(G,N))\rightarrow \cdots .$$
As $H_4(K(G,N))=H_3(G,N)$, $M(G,N)=H_3(K(G,N))$, $H_3(X)=M\cap
N/[M\cap N,G]$ and $H_2(K(G,N))=0$ (see $[2]$), we obtain the
following exact sequence
$$H_3(G,N)\rightarrow H_3(G/M,MN/M)\rightarrow H_4(X)\rightarrow
M(G,N)$$$$M(G/M,MN/M)\rightarrow M\cap N/[M\cap N,G][M,N]$$$$\ \ \ \
\ \ \ \ \ \ \ \ \ \rightarrow N/[N,G]\rightarrow
NM/M[N,G]\rightarrow 0.$$ Now, the Schur multiplier of the triple
$(G,M,N)$ is defined to be the fourth homology group of the space
$X$, $M(G,M,N)=H_4(X)$. If $K(G,M,N)$ denotes the mapping cone of
the canonical map $K(G,N)\rightarrow K(G/M,MN/M)$, then it is easy
to see that the space $X$ and $K(G,M,N)$ have the same homotopy
type. Therefore the above definition of the Schur multiplier of a
triple is coincide with the definition presented in $[4]$ by Ellis.
Using the above topological interpretation of the Schur multiplier
of a triple and similar to Theorem $2.7$, we can study the behavior
of the Schur multiplier of a triple with
respect to direct limits as follows.\\
\textbf{\bf Theorem 3.1.} \textit{Let $\{(G_i,M_i,N_i)\}_{i\in
I}$ be a directed system of triples of abelian groups such that
$M_i\cap N_i=1$ for all $i\in I$. Then}
$$\displaystyle{M(\lim_{\rightarrow} G_i, \lim_{\rightarrow} M_i,
\lim_{\rightarrow}N_i)\cong \lim_{\rightarrow} M(G_i,M_i,N_i).}$$
\textbf{Proof.} First, by the proof of Theorem 2.7, we recall
that for any abelian group $G_i$, the space
$\displaystyle{\lim_{\rightarrow} K(G_i,H_i)}$ can be considered
as
$\displaystyle{K(\lim_{\rightarrow}G_i,\lim_{\rightarrow}H_i)}$.
For the group $\displaystyle{\lim_{\rightarrow}M(G_i,M_i,N_i)}$,
we have the following diagram
\[\begin{array}{ccccccc}
 \ \ \ \ \ \displaystyle{\lim_{\rightarrow} K(G_i,N_i) \longrightarrow  \lim_{\rightarrow}K(\frac{G_i}{M_i},\frac{M_iN_i}{M_i}) } \\
 \hspace{0.1cm}  \downarrow \st{}{} \hspace{1.2cm} \downarrow
 \st{}{}\\
 \ \ \ \ \ \ \ 1 \longrightarrow \displaystyle{\lim_{\rightarrow}X_i.} \\
\end{array}\]
But by the above facts, we can replace this diagram by the new one
\[\begin{array}{ccccccc}
 \ \ \ \ \ \displaystyle{ K(\lim_{\rightarrow}G_i,\lim_{\rightarrow}N_i) \longrightarrow  K(\lim_{\rightarrow}\frac{G_i}{M_i},\lim_{\rightarrow}\frac{M_iN_i}{M_i}) } \\
 \hspace{0.1cm}  \downarrow \st{}{} \hspace{1.2cm} \downarrow
 \st{}{}\\
 \ \ \ \ \ \ \ 1 \longrightarrow \displaystyle{\lim_{\rightarrow}X_i.} \\
\end{array}\]
Now by the property of direct limit which preserves exact sequences,
we reform the above diagram to the following
\[\begin{array}{ccccccc}
 \ \ \ \ \ \displaystyle{ K(\lim_{\rightarrow}G_i,\lim_{\rightarrow}N_i) \longrightarrow  K(\frac{\displaystyle{\lim_{\rightarrow}G_i}}{\displaystyle{\lim_{\rightarrow}M_i}},
 \frac{\displaystyle{\lim_{\rightarrow}M_i\lim_{\rightarrow}N_i}}{\displaystyle{\lim_{\rightarrow}M_i}}) } \\
 \hspace{0.1cm}  \downarrow \st{}{} \hspace{1.2cm} \downarrow
 \st{}{}\\
 \ \ \ \ \ \ \ 1 \longrightarrow \displaystyle{\lim_{\rightarrow}X_i.} \\
\end{array}\]
The last diagram completes the proof. $\Box$\\

According to our topological definition of the Schur multiplier of a
pair of groups in previous section, it seems that we can generalize
this notion to triples more natural than the one of Ellis. In order
to define this new version of the Schur multiplier of a triple
$(G,M,N)$ consider the following homotopy pushout
\[\begin{array}{ccccccc}
 \ \ \ \ \ K(G) \longrightarrow  K(\frac{G}{N})  \\
 \hspace{0.4cm}  \downarrow \st{}{} \hspace{1cm} \downarrow
 \st{}{}\\
 \ \ \ \ K(\frac{G}{M}) \longrightarrow X. \\
\end{array}\]
Now we define the new version of the Schur multiplier of the
triple $(G,M,N)$ to be the third homology group of the space $X$,
$H_3(X)$.

Using the Mayer-Vietoris sequence for the above diagram we conclude
the following exact sequence $$\cdots\rightarrow
H_3(K(G))\rightarrow H_3(K(\frac{G}{N}))\oplus
H_3(K(\frac{G}{M}))\rightarrow H_3(X)\rightarrow
H_2(K(G))\rightarrow
$$$$H_2(K(\frac{G}{N}))\oplus H_2(K(\frac{G}{M}))\rightarrow
H_2(X)\rightarrow \cdots .$$ By Hopf Formula $[7]$ for any group
$G$, $H_2(K(G))=M(G)$; so we obtain the following exact sequence
$$\cdots\rightarrow H_3(K(G))\rightarrow H_3(K(\frac{G}{N}))\oplus
H_3(K(\frac{G}{M}))\rightarrow
$$$$M(G,N,M)\rightarrow M(G)\rightarrow M(\frac{G}{N})\oplus
M(\frac{G}{M})\rightarrow\cdots.$$
\textbf{\bf Remarks 3.2.}\\
$(i)$ Note that if $G=MN$, then by [2] $H_3(X)\cong ker(N\wedge
M\rightarrow G)$. Also, Ellis [4] mentioned that if $G=MN$, then
$M(G,N,M)\cong ker(N\wedge M\rightarrow G)$. Hence our definition of
the Schur multiplier of a triple of groups coincide with the
definition of Ellis, if $G=MN$.\\
$(ii)$ Since our definition of the Schur multiplier of a triple of
groups is a very natural generalization of the pair's one, we can
present more behaviors of this new notion with respect to free,
amalgamated free, and direct products and also a better behavior
with respect to direct limits than the one of Ellis.

The following results are evidences for the above claim.\\
\textbf{\bf Theorem 3.3.} \textit{For any two triple of groups
$(G_1,N_1,M_1)$ and $(G_2,N_2,M_2)$, we have the following
isomorphism}
$$M(G_1*G_2,
 \langle N_1*N_2\rangle^{G_1*G_2},\langle M_1*M_2\rangle^{G_1*G_2})\cong M(G_1,N_1,M_1)\oplus M(G_2,N_2,M_2).$$
\textbf{ Proof.} Suppose that $X_i$ $(i=1,2)$ is the pushout of
corresponding diagram
\[\begin{array}{ccccccc}
 \ \ \ \ \ K(G_i) \longrightarrow  K(\frac{G_i}{N_i})  \\
 \hspace{0.4cm}  \downarrow \st{}{} \hspace{1.2cm} \downarrow
 \st{}{}\\
 K(\frac{G_i}{M_i}) \longrightarrow X_i. \\
\end{array}\]
Similar to the proof of the Theorem $2.3$, using the fact that any
two direct limits commute, we conclude that the space
$X_1\vee X_2$ is also a pushout for the following diagram
\[\begin{array}{ccccccc}
\ \ \ \ \ \ \ \ \ \ \ K(G_1)\vee K(G_2) \longrightarrow  K(\frac{G_1}{N_1})\vee K(\frac{G_2}{N_2})  \\
 \hspace{0.6cm}  \downarrow \st{}{} \hspace{2.5cm} \downarrow
 \st{}{}\\
 \ K(\frac{G_1}{M_1})\vee K(\frac{G_2}{M_2}) \longrightarrow X_1\vee X_2. \\
\end{array}\]
According to Van-Kampen Theorem, we rewrite the above diagram as
follows
\[\begin{array}{ccccccc}
\ \ \ \ \ \ \ \ \ \ \ \ \ \ \ \ \ \ \ \ \ \ \ \ \ \ \ \ K(G_1*G_2) \longrightarrow  K(\frac{G_1}{N_1}*\frac{G_2}{N_2})=K(\frac{G_1*G_2}{\langle N_1*N_2\rangle^{G_1*G_2}})  \\
 \hspace{0.2cm}  \downarrow \st{}{} \hspace{2.1cm} \downarrow
 \st{}{}\\
  K(\frac{G_1*G_2}{\langle M_1*M_2\rangle^{G_1*G_2}})= K(\frac{G_1}{M_1}*\frac{G_2}{M_2}) \longrightarrow X_1\vee X_2 . \ \ \ \ \ \ \ \ \ \ \ \ \ \ \ \ \ \ \ \ \ \ \ \ \ \ \ \\
\end{array}\]
Now by the definition and the Mayer-Vietoris sequence for the above
diagram we conclude the following isomorphism $$M(G_1*G_2,
 \langle N_1*N_2\rangle^{G_1*G_2}, \langle M_1*M_2\rangle^{G_1*G_2})=H_3(X_1\vee X_2)\cong $$
 $$H_3(X_1)\oplus H_3(X_2)=M(G_1,N_1,M_1)\oplus
M(G_2,N_2,M_2).\ \Box$$

Also, by naturalness of our new version of the Schur multiplier of triple, we can give more results about the Schur multiplier of triples
of groups deduced with the proofs similar to those of pairs as follows.\\
 \textbf{\bf Theorem 3.4.} \textit{Let $(G_1,N_1,M_1)$, $(G_2,N_2,M_2)$ be two
triples of groups and $H\leq
N_1\cap N_2$, $H\leq M_1\cap M_2$. Then the Schur multiplier of the triple of
amalgamated free products $(G_1*_H G_2,N_1*_H
N_2,M_1*_H M_2)$ satisfies the
following exact sequence}
$$\cdots\rightarrow H_3(K(G_1*_H G_2))\rightarrow H_3(\frac{G_1}{N_1})\oplus
H_3(\frac{G_2}{N_2})\oplus H_3(\frac{G_1}{M_1})\oplus
H_3(\frac{G_2}{M_2})\rightarrow $$$$M(G_1*_H G_2,N_1*_H
N_2,M_1*_H M_2)\rightarrow M(G_1*_H G_2)\rightarrow
$$ $$M(\frac{G_1}{N_1})\oplus M(\frac{G_2}{N_2})\oplus
M(\frac{G_1}{M_1})\oplus M(\frac{G_2}{M_2})\rightarrow\cdots .$$
\textbf{\bf Theorem 3.5.} \textit{For any two triples of groups
$(G_1,N_1,M_1)$ and $(G_2,N_2,M_2)$, the Schur multiplier of the triple of direct
products $(G_1\times G_2,N_1\times N_2,M_1\times M_2)$ satisfies the
following exact sequence}
$$H_3(K(G_1))\oplus(M(G_1)\otimes G_2^{ab})\oplus (G_1^{ab}\otimes
M(G_2))\oplus H_3(K(G_2))\oplus
Tor(G_1^{ab},G_2^{ab})\rightarrow$$$$
H_3(K(\frac{G_1}{N_1}))\oplus(M(\frac{G_1}{N_1})\otimes
\frac{G_2^{ab}N_2}{N_2})\oplus (\frac{G_1^{ab}N_1}{N_1}\otimes
M(\frac{G_2}{N_2}))\oplus H_3(K(\frac{G_2}{N_2}))\oplus$$$$
Tor(\frac{G_1^{ab}N_1}{N_1},\frac{G_2^{ab}N_2}{N_2})\oplus
H_3(K(\frac{G_1}{M_1}))\oplus(M(\frac{G_1}{M_1})\otimes
\frac{G_2^{ab}M_2}{M_2})\oplus $$$$(\frac{G_1^{ab}M_1}{M_1}\otimes
M(\frac{G_2}{M_2}))\oplus H_3(K(\frac{G_2}{M_2}))\oplus
 Tor(\frac{G_1^{ab}M_1}{M_1},\frac{G_2^{ab}M_2}{M_2})\rightarrow$$$$ M(G_1\times G_2,N_1\times N_2,M_1\times
M_2)\rightarrow M(G_1\times G_2)\rightarrow
M(\frac{G_1}{N_1}\times\frac{G_2}{N_2})\oplus
M(\frac{G_1}{M_1}\times\frac{G_2}{M_2})\rightarrow\cdots .$$
\textbf{\bf Remarks 3.6.} Some special cases of the above exact sequence are as follows:\\
$(i)$ If $\frac{G_i}{N_i}$ and $\frac{G_i}{M_i}$ are finite cyclic
groups of orders $m_i$ and $l_i$ ($i=1,2$), respectively, then we
have the following exact sequence $$ \cdots \rightarrow {\bf
Z}_{m_1}\oplus {\bf Z}_{m_2}\oplus {\bf Z}_{d}\oplus {\bf
Z}_{l_1}\oplus {\bf Z}_{l_2}\oplus {\bf Z}_{c}\rightarrow
$$ $$M(G_1\times G_2,N_1\times N_2,M_1\times M_2)\rightarrow $$$$
M(G_1)\oplus M(G_2)\oplus (G_1\otimes G_2)\rightarrow {\bf
Z}_{d}\oplus {\bf Z}_{c}\rightarrow\cdots , $$  where $d=gcd(m_1,m_2)$
and $c=gcd(l_1,l_2)$.
 Moreover, if
$G_1$ and $G_2$ are also cyclic, then we get the following exact
sequence $$ \cdots \rightarrow {\bf Z}_{m_1}\oplus {\bf
Z}_{m_2}\oplus {\bf Z}_{d}\oplus {\bf Z}_{l_1}\oplus {\bf
Z}_{l_2}\oplus {\bf Z}_{c}\rightarrow$$$$ M(G_1\times G_2,N_1\times
N_2,M_1\times M_2)\rightarrow $$$$ G_1\otimes G_2\rightarrow {\bf
Z}_{d}\oplus {\bf Z}_{c}\rightarrow\cdots .$$
$(ii)$ If $gcd(|\frac{G_1}{N_1}|,|\frac{G_2}{N_2}|)=1$ and
$gcd(|\frac{G_1}{M_1}|,|\frac{G_2}{M_2}|)=1$, then the following
exact sequence holds $$ \cdots \rightarrow
H_3(\frac{G_1}{N_1})\oplus H_3(\frac{G_2}{N_2})\oplus
H_3(\frac{G_1}{M_1})\oplus H_3(\frac{G_2}{M_2})\rightarrow
$$$$M(G_1\times G_2,N_1\times N_2,M_1\times M_2)\rightarrow
M(G_1)\oplus M(G_2)\oplus (G_1\otimes G_2)\rightarrow$$$$
M(\frac{G_1}{N_1})\oplus M(\frac{G_2}{N_2})\oplus
M(\frac{G_1}{M_1})\oplus M(\frac{G_2}{M_2})\rightarrow\cdots .$$
Moreover, if $G_i$ is a finite group ($i=1,2$) such that
$gcd(|G_1|,|G_2|)=1$, then we conclude the following exact sequence
$$ \cdots \rightarrow H_3(\frac{G_1}{N_1})\oplus H_3(\frac{G_2}{N_2})\oplus
H_3(\frac{G_1}{M_1})\oplus H_3(\frac{G_2}{M_2}) \rightarrow
$$$$M(G_1\times G_2,N_1\times N_2,M_1\times M_2)\rightarrow M(G_1)\oplus
M(G_2)\rightarrow$$$$ M(\frac{G_1}{N_1})\oplus
M(\frac{G_2}{N_2})\oplus M(\frac{G_1}{M_1})\oplus
M(\frac{G_2}{M_2})\rightarrow\cdots .$$
$(iii)$ Finally, if $G_i$ is a finite cyclic group ($i=1,2$) such that
$gcd(|G_1|,|G_2|)=1$, then we have the following exact sequence
$$ \cdots \rightarrow {\bf Z}_{m_1}\oplus {\bf Z}_{m_2} \oplus {\bf Z}_{l_1}\oplus {\bf
Z}_{l_2} \rightarrow $$$$M(G_1\times G_2,N_1\times N_2,M_1\times
M_2)\rightarrow 0.$$\\

Finally, we can use Lemma $2.9$ and deduce a similar
result for triples, as follows.\\
\textbf{\bf Theorem 3.7.} \textit{The Schur multiplier of a triple
 of groups commutes with the direct limit of a directed system of
 triples of abelian groups.}

\end{document}